\documentclass{amsart}
\usepackage{amssymb,latexsym,amscd,hyperref}
\usepackage{pb-diagram}

\theoremstyle{plain}
\newtheorem{theorem}{Theorem}
\newtheorem{corollary}{Corollary}
\newtheorem{lemma}{Lemma}

\newtheorem{remark}{Remark}
\newtheorem{algorithm}{Algorithm}

\theoremstyle{definition}

\theoremstyle{remark}

\newcommand {\Q}{{\mathbb{Q}}}
\newcommand {\C}{{\mathbb{C}}}

\newcommand {\N}{{\mathbb{N}}}

\newcommand {\PP}{\mathbb{P}}

\newcommand {\K}{{\mathcal{K}}}

\newcommand {\OO}{{\mathcal{O}}}

\newcommand {\ida}{{\mathfrak{a}}}
\newcommand {\idp}{{\mathfrak{p}}}

\newcommand{\Norm}        {{\mathcal N}}

\newcommand{\Gal}      {\mathop{\rm {Gal}}}
\newcommand{\Rad}      {\mathop{\rm {Rad}}}

\newcommand{\Cl}      {\mathop{\rm {Cl}}\nolimits}

\newcommand{\PGL}        {{\mathop{\rm PGL}}}
\newcommand{\rk}        {{\mathop{\rm rk}}}

\newcommand{\eps}{\epsilon}
\begin{document}

\title{The number of $S_4$--fields with given discriminant}

\author{J\"urgen Kl\"uners}
\email{klueners@mathematik.uni-kassel.de}
\address{Universit\"at Kassel, Fachbereich Mathematik/Informatik,
Heinrich-Plett-Str. 40, 34132 Kassel, Germany.}


\begin{abstract} 
  We prove that the number of quartic $S_4$--extensions of the
  rationals of given discriminant $d$ is $O_\eps(d^{1/2+\eps})$ for
  all $\eps>0$. For a prime number $p$ we derive that the dimension of
  the space of octahedral modular forms of weight 1 and conductor $p$
  or $p^2$ is bounded above by $O(p^{1/2}\log(p)^2)$.
\end{abstract}

\maketitle
 
\renewcommand{\theenumi}{(\roman{enumi})}

\section{Introduction}
For a number field $k$ we denote by $d_k\in\N$ the absolute value of
the field discriminant of $k$. The class group will be denoted by
$\Cl_k$ and the $p$-rank $\rk_p(A)$ of an abelian group $A$ is defined
to be the minimal number of generators of $A/A^p$.  We denote by
$\Norm$ the absolute norm. The symbol $O_\eps$ denotes the usual Landau
symbol $O$, where the implied constant is depending on $\eps$.

In this note we answer a question of Akshay Venkatesh about the number
of $S_4$--extensions of degree 4 with given discriminant $d$. It is
conjectured that this number is $O_\eps(d^\eps)$ for all $\eps>0$. In
average we have the stronger result (see \cite{Bha,Bel}):
$$\lim_{x\rightarrow\infty}\frac{1}{x}\sum_{K: d_K\leq x}1 = c(S_4),$$
where $K$ runs through all quartic $S_4$--extensions and $c(S_4)>0$ is
explicitly given.  We prove the bound $O_\eps(d^{1/2+\eps})$ for all
$\eps>0$ which improves the bound $O_\eps(d^{4/5+\eps})$ given in
\cite{MiVe}. 

As an application we give an upper bound for the dimension of the
space of octahedral forms of weight 1 and given conductor $N$. In the
general case the best known bound is $O_\eps(N^{4/5+\eps})$ for all
$\eps>0$ given in \cite{MiVe}. For squarefree conductors this bound is
improved to $O_\eps(N^{2/3+\eps})$ on average. In this note we are
able to prove the upper bound $O_{\eps}(N^{1/2+\eps})$ in many cases,
e.g. when $N$ is prime or a square.

The discrepancy between the expected bound $O_\eps(d^\eps)$ and the
proven bound $O_\eps(d^{1/2+\eps})$ for the number of
$S_4$--extensions of discriminant $d$ comes from the fact that we can
only use weak bounds for the $3$--rank of the classgroup of quadratic
fields and the $2$-rank of the classgroup of non-cyclic cubic fields.

In order to understand the problems which arise we give the following
easy example. Let us count the number of cubic $S_3$--extensions
$M/\Q$ of discriminant $d$ such that the normal closure contains a
given quadratic extension $k$. Since every unramified cyclic cubic
extensions $N/k$ corresponds to a cubic extension $M$ we see that the
number of elements $h_3$ of order $3$ in the classgroup $\Cl_k$ plays
an important role.  In the general case we can only use the estimate
$h_3\leq \#\Cl_k$ and the latter one can be bounded by
$O(d^{1/2}\log(d))$ using Lemma \ref{boundcl}.  It is very difficult
to improve this trivial bound for elements of order $p$ in the class
group when $p>3$. Just recently for $p=3$ Helfgott and Venkatesh
\cite{HeVe} ($\lambda=0,44179$) and independently Pierce \cite{Pie}
($\lambda=0.49108$ or $\lambda=0.41667$ in special cases) proved that
for all $\eps>0$ we get: $$3^{\rk_3(\Cl_k)} =
O_\eps(d_k^{\lambda+\eps}).$$
Using this improved bound it is
straightforward to get the upper bound $O_\eps(d^{\lambda+\eps})$ for
the number of cubic $S_3$--extensions.

In the following we would like to explain the idea of the proof of our
main result.  We will improve the following elementary approach given
in \cite[p.  101]{Duke}.  In the worst case we cannot exclude the case
that there exists a quadratic field $k/\Q$ such that
$3^{\rk_3(\Cl_k)}=O(d_k^{1/2}\log(d_k))$.  Using these unramified
$C_3$--extensions of $k$ there are $\frac{3^{\rk_3(\Cl_k)}-1}{3-1}$
non-cyclic cubic fields $M$ of the same discriminant. In the worst
case all these extensions have a large $2$--rank, i.e.
$2^{\rk_2(\Cl_M)}=O(d_M^{1/2}\log(d_M)^2)$.  Every unramified
$C_2$-extension leads to an $S_4$-extension $K$ of degree 4 of the
same discriminant $d_K=d_k=d_M$.  Using this idea we get the upper
bound $O(d_K\log(d_K)^3)$ for the number of $S_4$--extensions of
discriminant $d_K$.

As we see from the above example it is a problem for our upper estimates when
$\rk_2(\Cl_M)$ and $\rk_3(\Cl_k)$ are big. We will use Theorem \ref{up3}
proved by Frank Gerth III which says that $\rk_3(\Cl_M)$ has about the
same size as $\rk_3(\Cl_k)$. This means that $\rk_3(\Cl_M)$ is big when 
$\rk_3(\Cl_k)$ is big. This implies that $\rk_2(\Cl_M)$ must be small.

E.g. consider the special case that $d$ is squarefree, i.e. the
corresponding $S_3$--extension $L$ is unramified over $k$. Then the
first part of Theorem \ref{up3} and the above explained elementary
approach already proves our wanted result, i.e. the number of
$S_4$--extensions of disriminant $d$ is bounded by $O(d^{1/2+\eps})$.

\section{Parameterizing $S_4$--extensions}
\label{sec:main}
 Let $K/\Q$ be a
quartic field such that the normal closure $N$ has Galois group $S_4$.
Then there is a unique normal subfield $L$ of degree 6 having Galois
group $S_3$. We denote by $M$ a subfield  of $L$ of degree 3 and by $k$
the unique subfield of degree 2 of $L$ (or $N$).
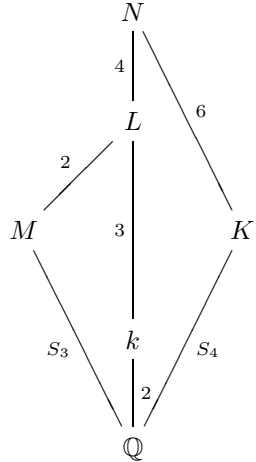
\[
\begin{diagram}
\node[2]{N} \arrow{s,l,-}{4}\arrow{sse,l,-}{6}\\
\node[2]{L} \arrow{sw,l,-}{2} \arrow[2]{s,l,-}{3}\\
\node{M} \arrow{sse,r,-}{S_3} \node[2]{K} \arrow{ssw,r,-}{S_4}\\
\node[2]{k} \arrow[1]{s,r,-}{2}\\
\node[2]{\Q}
\end{diagram}
\]
For $n\in\N$ we define $\Rad(n):=\prod\limits_{p \mid n}p$, 
where the product is only taken over primes.
  To each 
$K/\Q$ as above we associate a triple
$$(a,b,c)=(\Rad(d_k),\Rad(\Norm(d_{L/k})),\Rad(\Norm(d_{N/L})))\in\N^3$$
of squarefree numbers. We define
\begin{equation}
  \label{eq:fib}
  \Psi: \K \rightarrow \N^3, K \mapsto (a,b,c),
\end{equation}
where $\K$ is the set of quartic $S_4$--extensions of $\Q$ up to isomorphy.
$\Psi$ is a well defined mapping with bounded fibers. In the rest of this section
we want to give upper bounds for the size of the fibers, i.e. to give an upper bound
for the numbers of fields $K$ which are associated to a given triple $(a,b,c)$? 

Assuming this situation $k$ is one of the following quadratic fields.
If $2\nmid a$ we get that $k=\Q(\sqrt{\pm a})$ where the sign is
positive if $a\equiv 1 \bmod 4$.  If $2\mid a$ then $k$ is one of the
following three fields: $\Q(\sqrt{a}), \Q(\sqrt{-a})$, and
$\Q(\sqrt{\pm a/2})$, where the sign is positive when $a/2 \equiv 3
\bmod 4$.  Therefore at most 3 quadratic fields are associated to a
given $a$.  The number of $b$'s for a given field $k$ can be easily
bounded by the following lemma. In the following we denote by
$\omega(b)$ the number of prime factors of $b$.
\begin{lemma}\label{lem:rank3}
  Let $b\in\N$ as above. Then all fields $M$ (up to isomorphism) such
  that $L/K$ is only ramified in primes dividing $b$ are contained
  in the ray class field of $\ida:=3b\OO_k$. The number of those
  extensions can be bounded by
$$\frac{3^r-1}{3-1}, \mbox{ where }r=\rk_3(\Cl_k)+\omega(b)+2.$$
\end{lemma} 
\begin{proof}
  We are looking for all fields which are at most ramified in primes
  dividing $b$. We need to choose $\ida$ in such a way that all these
  fields are subfields of the ray class field of $\ida$.
  For primes $\idp$ not dividing 3 it is sufficient that
  $\idp\mid \ida$. For the wildly ramified primes there exists a
  maximal exponent such that all these fields occur as subfields
  \cite[p. 58]{Serre} of the ray class field of $\ida$.  Using
  elementary properties of the ray class group $\Cl_\ida$ we get that
  $$\rk_3(\Cl_\ida) \leq \rk_3(\Cl_k) + \rk_3((\OO_k/\ida)^*).$$
  For all prime ideals $\idp$ not dividing 3 we get that the $3$-rank of
  $(\OO_k/\idp)^*$ is at most 1 which shows that
  $\rk_3(\OO_k/p\OO_k)\leq 2$.  Equality can only occur in the case
  $p\equiv 1 \bmod 3$, where $p\in\PP\cap \idp.$   
  In this case there exists a $C_3$--extension of
  $\Q$ only ramified in $p$. Denote by $A$ the $3$--part of the ray class group
  $\Cl_{\ida}$. We can write
  $A:=A^+\oplus A^-$, where the classes in $A^+$ are invariant under $\Gal(k/\Q)$.
  Because a prime $p\equiv 1 \bmod 3$ increases the $3$--rank of $A^+$ by one, we
  get that all odd primes increase the $3$--rank of $A^-$ by at most one.
  The theory used in
  \cite[Section 6]{FiKl} shows that $S_3$--extensions correspond to quotients
  of index 3 of $A^-$.
  Finally we
  need to estimate the 3--rank for $\OO_k/\idp^w$ for primes dividing
  $3$.  In \cite{HePaPo} it is proved that the $p$-rank of
  $(\OO_k/\idp^w)^*$ is at most $[k_\idp:\Q_p]+1$. In all cases it is
  sufficient to add 2 since there is one $C_3$-extension of $\Q$ only
  ramified in $3$.
\end{proof}
We use the trivial class group bound which can be found
in \cite[Theorem 4.4]{Nar}.
\begin{lemma}\label{boundcl}
  For all $n\in\N$ there exists a constant $c(n)$
  such that for all number fields $F$ of degree $n$ we have:
  $$|\Cl_F| \leq c(n) d_F^{1/2}\log(d_F)^{n-1}.$$
\end{lemma}
Trivially, we have $3^{\rk_3(\Cl_k)}\leq |\Cl_k|$.  For
a given cubic $S_3$-field $M$ we prove a similar lemma as Lemma \ref{lem:rank3}.
\begin{lemma}
  Let $c\in\N$ be as above. Then
  the number of $S_4$-extensions $N$ which contain a given $S_3$-field
  $M$ such that $\Norm(d_{N/L})$ is only divisible by primes dividing
  $c$ is bounded by $$2^r-1, \mbox{ where
  }r=\rk_2(\Cl_M)+3\omega(c)+6.$$
\end{lemma}
\begin{proof}
  In \cite[Lemmata 4,5]{Bai} it is proven that the Galois closure of
  $M(\sqrt{\alpha})$ for $\alpha\in M$ has Galois group $S_4$ if and only if
  $\Norm(\alpha)$ is a square. If $\Norm(\alpha)$ is a square this
  certainly implies that the norm of the principal ideal $(\alpha)$ 
  is a square. Therefore we get an upper bound if we count all extensions
  such that the conductor is a square. For a prime $p\ne 3$ we have
  at most three possibilities to produce squarefree ideals of norm $p^2$.
  The $6$ is computed in a similar way as in Lemma \ref{lem:rank3} and
  gives an upper bound for the contribution of primes above 3.
\end{proof}
Altogether we get the following upper bound for the number of
$S_4$-fields associated to a given triple $(a,b,c)$: 
\begin{equation}
  \label{eq:number}
3(\frac{3^{r_1}-1}{3-1})(2^{r_2}-1)
\leq 3/2\cdot 9 \cdot 2^6 3^{\rk_3(\Cl_k)} 2^{\rk_3(\Cl_M)}
3^{\omega(b)} 8^{\omega(c)},  
\end{equation}
where $r_1=\rk_3(\Cl_k)+\omega(b)+2,\;\;r_2=\rk_3(\Cl_M)+3\omega(c)+6.$

The following theorem relates the $3$--parts of the classgroups of $k$ and $M$.
\begin{theorem}(Gerth III)\label{up3}
  Let $M/\Q$ be a non-cyclic cubic extension and denote by $L$ the
  normal closure of $M$ and by $k$ the unique quadratic subfield of
  $L$. Then the following holds.
  \begin{enumerate}
  \item If $L/k$ is unramified, then $\rk_3(\Cl_M)=\rk_3(\Cl_k)-1$.
  \item $\rk_3(\Cl_M) = \rk_3(\Cl_k)+ t -1 -z - y$, where $y\leq t-1$ and $t$ is the number of prime ideals of $\OO_k$
    which ramify in $L$. Furthermore we have $0\leq z\leq u$ where $u$ is the number of primes which are totally ramified in
    $M$ but split in $k$.
  \item $\rk_3(\Cl_M) \geq  \rk_3(\Cl_k) -u$
  \end{enumerate}
\end{theorem}
\begin{proof}
  The first part is Theorem 3.4 in  \cite{Ger}. The second part is Theorem 3.5. The last part is an immediate consequence.
\end{proof}

Since we are only interested in the asymptotic behaviour we can ignore ramification in $2$ and
$3$. Therefore we define $S:=\{2,3\}$ and 
$a^S$ to be the largest number dividing $a$ which is
coprime to $S$. Using this we easily see that $d_M^S =a^S(b^S)^2$, where $M$ is one of the cubic
extensions constructed above.
Using Theorem \ref{up3} we get the
following estimate for $3^{\rk_3(\Cl_k)} 2^{\rk_3(\Cl_M)}$.
\begin{lemma}\label{rankrel}
  Let $M,k$ be the fields defined before.
  Then there exists a constant $C>0$ such that
  $$3^{\rk_3(\Cl_k)} 2^{\rk_2(\Cl_M)} \leq C a^{1/2}b \log(ab^2)^2
  3^{\omega(b)}.$$
\end{lemma}
\begin{proof}
  Theorem \ref{up3} shows $\rk_3(\Cl_M)\geq \rk_3(\Cl_k)-\omega(b)$. Therefore we get:
  $$3^{\rk_3(\Cl_k)}2^{\rk_2(\Cl_M)} \leq
  3^{\rk_3(\Cl_M)}3^{\omega(b)}2^{\rk_2(\Cl_M)}\leq
  3^{\omega(b)}|\Cl_M|.$$
  Using Lemma \ref{boundcl} and the fact that
  $d_M^S=(ab^2)^S$ differs from $d_M$ by something which can be
  bounded by a constant we get the desired bound.
\end{proof}
Combining Lemma \ref{rankrel} and \eqref{eq:number} we deduce the following corollary.
\begin{corollary}\label{cor:upp}
    The number of elements of the fiber $\Psi^{-1}(a,b,c)$ is bounded by
    $$3^32^5 C a^{1/2}b\log(ab^2)^2 9^{\omega(b)}8^{\omega(c)}.$$
\end{corollary}

\section{Upper bounds for quartic $S_4$--extensions with given discriminant}
\label{sec:upp}
In this section we prove an upper bound for the number of quartic
$S_4$--extensions with given discriminant. In order to do this we need
to compute the discriminant $d_K$ using the triple $(a,b,c)$.  In a
second step we determine how many triples may lead to the same
discriminant.

Let us assume that we have given a field $K\in\K$ with $\Psi(K)=(a,b,c)$ ramified in $p$.
Assuming $p\ne 2,3$ we can compute the cycle
shape of a generator of the cyclic inertia group at $p$ in the degree 4
representation of $S_4$.   Here we denote by the cycle shape the
length of the cycles if we decompose a group element into disjoint
cycles.  Using local theory we get for primes $p>3$ the following
identities, where $v_p$ denotes the ordinary $p$-valuation.  The
results are given in the following table:
\begin{center}
\begin{tabular}[c]{|l|l|l|}
\hline
& cycle shape & $v_p(d_K)$ \\ \hline
$p\mid a, p\nmid bc$ & $1^2 2$ & 1 \\ \hline
$p\mid a, p \mid c, p\nmid b$ & $4$ & 3 \\ \hline
$p\mid b, p\nmid ac$ & $1 3$ & 2 \\ \hline
$p\mid c, p\nmid ab$ & $2^2$ & 2 \\ \hline
\end{tabular}
\end{center}
The other cases cannot occur since in these cases the inertia group
would not be cyclic. The cases $p=2$ or $p=3$ can be handled by
analyzing the local Galois groups. We still use the definition $a^S$ for $S:=\{2,3\}$ from
the preceding section and get:
$$d_K^S = a^S(b^S)^2(c^S)^2.$$
The contribution of the primes 2 and 3
is bounded by a constant factor. Therefore we ignore these primes in
the following.

Using the results of the preceding section it remains to 
count the number of triples $(a,b,c)$ which may lead to the same discriminant.
In the following let $d$ be a discriminant of a quartic
$S_4$--extension, 
\begin{theorem}
  Let $d=2^{e_2}3^{e_3}d_1d_2^2d_3^3$ such that $6d_1d_2d_3$ is squarefree.
  Then the number of $S_4$-fields with discriminant $d$ is bounded above by 
  \begin{enumerate}
   \item $\tilde{C} (d_1d_3)^{1/2}d_2\log(d_1d_3d_2^2)^2 18^{\omega(d_2)}8^{\omega(d_3)}$
 for a suitable $\tilde C>0$.
   \item $O_\eps(d^{1/2+\eps})$ for all $\eps>0$.
  \end{enumerate}
\end{theorem}
\begin{proof}
  Using the above discussion all fields $K/\Q$ with $\Psi(K) = (a,b,c)$ have the property:
  $$a^S=d_1d_3, \;d_3 \mid c^S \mbox{ and }(bc)^S = d_2d_3.$$
  Therefore we have $2^{\omega(d_2)}$ possibilities for choosing $b^S$. 
  The number of possibilities for the $2$ and the $3$--part can be bounded by a constant.
  Using Corollary \ref{cor:upp} the worst case is when $b^S=d_2$ and therefore we get for some computable
  constant $\tilde{C}>0$
  $$\tilde{C} 2^{\omega(d_2)}  (d_1d_3)^{1/2}d_2\log(d_1d_3d_2^2)^2 9^{\omega(d_2)}8^{\omega(d_3)}$$
  as an upper bound. For the second statement we write
 $x^{\omega(d)}=O(d^\eps)$ for a given number $x$ and get the desired result.
\end{proof}
\begin{remark}
  For squarefree discriminants $d$ we can derive  the better upper bound $O(d^{1/2}\log(d)^2)$.
\end{remark}

We can combine this result with well known results to get bounds for degree 4 fields.
\begin{theorem}
  The number of degree 4 fields of given discriminant $d$ is bounded
  above by $O_\eps(d^{1/2+\eps})$ for all $\eps>0$.
\end{theorem}
\begin{proof}
  Using the theorem of Kronecker-Weber we easily get that the number
  of fields with Abelian Galois group is bounded by $O_\eps(d^\eps)$
  for all $\eps>0$. Since $D_4$-fields can be constructed by quadratic
  extensions over quadratic extensions and the $2$-torsion part of the
  class group can be easily controlled, we get the same result for
  $D_4$-extensions. For $A_4$-extensions we use the same approach as
  in the $S_4$-case. The main difference is that we have only one step
  where we have to consider class groups. This gives
  $O_\eps(d^{1/2+\eps})$ for the number of such extensions with given
  discriminant $d$. Using more advanced methods \cite{MiVe} this
  number can be reduced to $O_\eps(d^{1/3+\eps})$.
\end{proof}

\section{Upper bounds for the dimension of the space of octahedral modular forms of given conductor}
\label{sec:mod}

In this section we give upper bounds for the dimension of the space of octahedral
modular forms of weight 1. Denote by $G_\Q$ the absolute galois group $\Gal(\bar\Q/\Q)$.
Suppose we have given a quartic $S_4$--extension $K/\Q$ which gives rise
to a projective representation $\tilde\rho: G_\Q \rightarrow \PGL_2(\C)$. The conductor
of this projective representation is defined to be the product of the local conductors
of $\rho_{|G_{\bar{\Q}_p}}: G_{\bar{Q}_p} \rightarrow \PGL_2(\C)$ which is the minimal $p$--power
of a so--called local lift, see e.g. \cite[\S 6]{Ser77} or \cite{Wong} for more details.

In this section we count $S_4$--extensions using the above defined conductor. A
prime $p$ divides the conductor if and only if $p$ divides the
discriminant. To simplify all computations we ignore the contribution
of the 2-- and 3--part of the conductor. For all other (tamely)
ramified primes we have the property that $p$ exactly divides the
conductor when the local Galois group is cyclic. Otherwise the local 
Galois group is dihedral and we get that $p^2$ exactly divides
the conductor \cite[Prop. 1, p. 144]{Wong}.

To each projective representation with image $S_4$ we can associate an octahedral modular form of
the same conductor. This means that we get the corresponding bounds for the modular forms
when we compute the bounds for the number of projective representations 
(see e.g. \cite{Duke,Wong} for more details).

In order to use the results of Section \ref{sec:main} we need to compute
the conductor of the associated modular form only using the triple $(a,b,c)$. Similar to the
discriminant case we can do all computations locally. In the discriminant case it was only important
to know the inertia group. Now it is important to know the decomposition group. Let $p>3$ be a divisor
of $abc$. Then $p$ exactly divides the conductor if the decomposition group is cyclic.  
In the following table we collect the information we 
get (for $p>3$) using the prime ideal factorization $p\OO_K=\prod_{i=1}^r \idp_i^{e_i}$.
We remark that some cases can be distinguished by congruence conditions. In  the last column
we denote the letters which are divisible by $p$. The information $v_p(d)$ on the discriminant is not needed
in this section.
\begin{center}
\begin{tabular}{|l|l|l|l|l|l|l|}
  \hline
  & $D_p$ & $I_p$ & $v_p(N)$ & $v_p(d)$  & & $p\mid$ \\ \hline 
$\idp_1^2 \idp_2\idp_3$  & $C_2$ & $C_2$ & 1 & 1  & & $a$ \\ \hline 
$\idp_1^2 \idp_2$ & $C_2\times C_2$ & $C_2$ & 2 & 1 & & $a$ \\ \hline
$\idp_1^2$ & $C_2\times C_2$ or $C_4$ & $C_2$ & 2 or 1 & 2  & & $c$ \\ \hline
$\idp_1^2\idp_2^2$ & $C_2\times C_2$ or $C_2$ & $C_2$ & 2 or 1 & 2  & & $c$ \\ \hline 
$\idp_1^4$ & $D_4$ & $C_4$ & 2 & 3 & $p\equiv 3 \bmod 4$ & $a,c$ \\ \hline
$\idp_1^4$ & $C_4$ & $C_4$ & 1 & 3 & $p\equiv 1 \bmod 4$ & $a,c$ \\ \hline
$\idp_1^3\idp_2$ & $C_3$ & $C_3$ & 1 & 2& $p\equiv 1 \bmod 3$ & $b$ \\ \hline
$\idp_1^3\idp_2$ & $D_3$ & $C_3$ & 2 & 2& $p\equiv 2 \bmod 3$ & $b$ \\ \hline
\end{tabular}
\end{center}
Let $K/\Q$ be a quartic $S_4$--extension with associated triple $(a,b,c)$ and conductor $N$. Then
we write $a=a_1a_2$, $b=b_1 b_2$, $c=c_0c_1c_2$, where $c_0:=\gcd(a,c)$ such that
$$N^S = (a_1a_2^2 b_1b_2^2 c_1 c_2^2)^S.$$
Since $\gcd(b,ac)^S=1$ we easily see that
$a_1^S,a_2^S,b_1^S,b_2^S,c_1^S,c_2^S$ are 
pairwise coprime.
Using the above table we know that $b_i$ is (up to the $3$--part) exactly divisible by the primes
dividing $b$ which are congruent to $i \bmod 3 \;(i=1,2)$.
\begin{theorem}
  Let $N=2^{n_2}3^n_3 N_{1,1}N_{1,2}N_2^2$ such that $6N_{1,1}N_{1,2}N_2$ is squarefree.
  Furthermore we assume that $p\mid N_{1,i}$ if and only if $p\equiv i \mod 3$ ($i=1,2$).
  Then the number of $S_4$--fields of given conductor $N$ is bounded above by
  $$C 54^{\omega(N)} N_{1,1}N_{1,2}^{1/2}N_2\log(N)^2$$ 
  for a suitable $C>0$.
\end{theorem}
\begin{proof}
  We have $3^{\omega(N)}$ possibilities to partition the primes into three sets corresponding
  to $a,b,c$. Furthermore we have at most $2^{\omega(N)}$ possibilities for $c_0$.
  Using Corollary \ref{cor:upp} we have the worst case when $b$ is big. Primes dividing
  $N_{1,2}$ cannot divide $b$. Here we get the worst case when these primes divide $a$.
  Therefore we get as an upper bound:
  $$\tilde{C} 3^{\omega(N)}2^{\omega(N)}  N_{1,1}N_{1,2}^{1/2}N_2 \log(N_{1,2}N_{1,1}^2N_2^2)^2 9^{\omega(N)}.$$
  We easily get the desired result.
\end{proof}

To get good estimates for the dimension of the space of octahedral forms with given
conductor we have to avoid that $b_1$ is big.  Using this we can
derive the following corollaries.
\begin{corollary}
  Let $p$ be a prime. Then the dimension of the space of octahedral modular forms of weight 1 and conductor
  $p$ or $p^2$ is bounded above by $O(p^{1/2}\log(p)^2)$.
\end{corollary}
\begin{proof}
  The quadratic subextension must be ramified in at least one prime. Therefore $p\mid a$ for all
  possible triples.
\end{proof}
\begin{corollary}
  Assume that all primes which exactly divide
  $N$ are congruent to $2 \bmod 3$.  Then the dimension of the space of octahedral
  forms of weight 1 and conductor $N$ is bounded above by $O(N^{1/2+\eps})$
  for all $\eps>0$.
\end{corollary}
\begin{proof}
  We have $N_{1,1}=1$ and the assertion follows.
\end{proof}

This improves the bound $O(N^{4/5+\eps})$ given in \cite{MiVe}. We remark that in the case that
$b_1$ resp. $N_{1,1}$ is big we only get the trivial linear bound using our method.

\section*{Acknowledgments}
I thank Karim Belabas and Gunter Malle for fruitful discussions and reading
a preliminary version of the paper.
\bibliographystyle{amsalpha} 
\bibliography{myref} 
\end{document}